\documentclass[twocolumn,letterpaper]{IEEEAerospaceCLS}  % only supports two-column, letterpaper format

\usepackage{amssymb}
\usepackage{amsmath}
\usepackage{mathtools}
\usepackage[]{graphicx}    % We use this package in this document
\newcommand{\ignore}[1]{}  % {} empty inside = %% comment
\def\C{\boldsymbol{C}}

\begin{document}
\title{Nonlinear Optimal Guidance for Cooperatively Imposing Relative Intercept Angles}

\author{%
Han Wang\\ 
School of Aeronautics and Astronautics\\
Zhejiang University\\
Hangzhou, China, 310027\\
han\_wang98@zju.edu.cn
\and 
Zheng Chen\\
School of Aeronautics and Astronautics\\
Zhejiang University\\
Hangzhou, China, 310027\\
Huanjiang Laboratory\\
Zhuji, Zhejiang, China, 311899\\
z-chen@zju.edu.cn
%%%% IMPORTANT: Use the correct copyright information--. %%%%%
\thanks{ }           % This creates the copyright info that is the correct 2024 data.
%\thanks{{}}         % Use this copyright notice only if you are employed by ...
%\thanks{{}}          % Use this copyright notice only if you are employed by ...
%\thanks{{}}    % Use this copyright notice is you are employed by ...
}

\maketitle

\thispagestyle{plain}
\pagestyle{plain}

\maketitle

\thispagestyle{plain}
\pagestyle{plain}

\begin{abstract}
The optimal cooperative guidance in the nonlinear setting for intercepting a target by multiple pursuers is studied in the paper. As certain relative angles can improve  observability, the guidance command is required to cooperatively control the pursuers to intercept the target with specific relative angles. By using the neural networks, an approach for real-time generation of the nonlinear cooperative optimal guidance command is developed. Specifically, the optimal control problem with constraints on relative intercepting angles is formulated. Then, Pontryagin's maximum principle is used to derive the necessary conditions for optimality, which are further employed to parameterize the nonlinear optimal guidance law. As a result, the dataset for the mapping from state to nonlinear optimal guidance command can be generated by a simple propagation. A simple feedforward neural network is trained by the dataset to generate the nonlinear optimal guidance command. Finally, numerical examples are presented, showing that a nonlinear optimal guidance command with specific relative angles can be generated within a faction of a millisecond. 
\end{abstract}

\tableofcontents

%%%%%%%%%%%%%%%%%%%%%%%%%%%%%%%%%%%%%%
\section{Introduction}
%%%%%%%%%%%%%%%%%%%%%%%%%%%%%%%%%%%%%%
%In complex environment, one may encounter difficulties, such as close-in weapon system (CIWS) equipped on surface ship, when intercepting a high-value target by a single pursuer. 
The research on intercepting scenarios with multiple pursuers is becoming increasingly popular in recent years. Because single pursuer cannot satisfy the demand of various application, a cooperative guidance architecture for multiple pursuers manifests the superiority in temporal and spatial scale. Thus, a temporal cooperative scheme involving multiple pursuers is considered as an effective solution and has drawn much attention due to its improved interception performance and its high level of fault tolerance \cite{jeon2006impact}. 

However, pure temporal cooperation for simultaneous
arrival can not release the full effectiveness of multiple pursuers. Therefore, spatial cooperation for encirclement is also investigated  \cite{shaferman2021near}-\cite{zhai2016coverage}. Accordingly, a spatial cooperation for multi-directional interception is devised to  effectively enhance the observability and limit the evasion of the opponent. Thus, the spatial cooperative attack helps to enhance the target identification and decrease its evasive probability. 

The spatial cooperative interception enables the group of pursuers to formulate superior geometry relative to the target. Eventually, the purpose of spatial cooperation is to control the different angles of the pursuer-target collision triangle. In \cite{shaferman2015cooperative}, an explicit cooperative guidance law by controlling the relative intercept angles to formulate angular geometry based on the linearized kinematics was proposed. And it is proved that explicit cooperation yields better performance than implicit ones. Fonod and Shima developed a cooperative guidance law to improve the observability of the target by imposing desired relative flight-path angles \cite{fonod2017estimation}. Meanwhile, in the game of pursuer and target, the evasive maneuver of the target may result in a large heading angle. Under these engagement conditions, the small angle assumption is violated and the performance of aforementioned algorithms that are based on linearized kinematics will deteriorate dramatically.

 The nonlinear control theory has been used to develop cooperative guidance strategy to avoid the issue of linear optimal cooperative guidance laws. A sliding-mode control theory was leveraged by Song et al. \cite{song2017three} to realize implicit coordination of the line-of-sight (LOS) angle. Finite-time consensus theory and the super-twisting control algorithm were employed by Zhang et al. \cite{zhang2020finite} to control each pursuer to reach the target with the desired LOS angle. In \cite{li2022optimal}, the concept of encirclement was introduced in the field of cooperative guidance, and a pseudocontrol-effort optimal encirclement guidance law was derived in a relative virtual frame. Nonetheless, to the author's best knowledge, it is not common to see the studies on real optimal explicit cooperative guidance laws.

In this paper, based on the work in \cite{chen2019nonlinear,doi:10.2514/1.G006666}, the focus is on developing an explicit cooperative guidance law that is able to guide each pursuer to reach its desirable relative intercept angle with minimum control effort. By parameterizing the nonlinear optimal guidance law using Pontryagin’s Maximum Principle (PMP) \cite{pontryagin2018mathematical}, we can generate the dataset for the mapping from state to nonlinear optimal guidance command by solving some simple initial value problems. After training a basic feedforward neural network (FNN) with the dataset, it can swiftly produce an optimal guidance command within a fraction of a millisecond.

%%%%%%%%%%%%%%%%%%%%%%%%%%%%%%%%%%%%%%%%%%%
\section{Problem Formulation}
%%%%%%%%%%%%%%%%%%%%%%%%%%%%%%%%%%%%%%%%%%%
\subsection{Engagement Geometry}
\begin{figure}\label{OneColumn}
	\centering
	\includegraphics[width=3.25in]{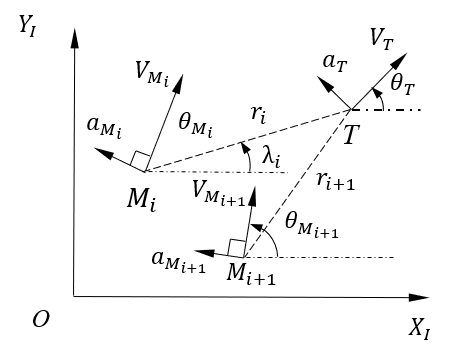}\\
	\caption{\textbf{ Engagement geometry.}}
	\label{Geometry}
\end{figure}

The planar engagement geometry for $N$ pursuers against one maneuvering target is depicted in Fig.~\ref{Geometry}. The inertial frame is denoted by $X_IOY_I$, $M_i$ signifies the $i$th pursuer, and $T$ represents the target. During the terminal guidance stage, the engaged pursuers are assumed to have a constant velocity. The speed is denoted by $V$, and the acceleration is represented by $a$. Symbol $\theta$ is called the heading angle. Those notations with subscript $M_i$ and $T$ refer to corresponding values for the $i$th pursuer and the target, respectively.

The terminal intercept angles of the $i$th pursuer and the $(i+1)$th pursuer are $\theta_{M_i}+\theta_T$ and $\theta_{M_{i+1}}+\theta_T$, respectively. The difference between these two intercept angles is the relative intercept angle between the two pursuers from the target’s perspective. This is the angle that will be forced between two consecutive pursuers. Note that, if the target does not maneuver between the two intercepts or the intercepts occur at the same time, this angle will be $\theta_{M_{i+1}}-\theta_{M_i}$ \cite{shaferman2015cooperative}.

The differential equations for the relative motion are given by
\begin{equation}
	\begin{split}
		\dot r_i &= V_T\cos(\theta_T-\lambda_i)-V_{M_i}\cos(\theta_{M_i}-\lambda_i)\\
		\dot \sigma_i &= \frac{1}{r_i}\big[V_T\sin(\theta_T-\lambda_i)-V_{M_i}\sin(\theta_{M_i}-\lambda_i)\big]\\
		\dot{\theta}_{M_i} &= \frac{a_{M_i}}{V_{M_i}}\\
		\dot{\theta}_{T} &= \frac{a_T}{V_{M_i}}
	\end{split}
\end{equation}
where $r_i$ represents the range between the $i$th pursuer and the target, $\lambda_i$ denotes the LOS angle.
 
The nonlinear kinematics for the $i$th pursuer is expressed as
\begin{align}
	\begin{split}
		\dot x_i(t)&=\cos\theta_{M_i}(t)\\
		\dot y_i(t)&=\sin\theta_{M_i}(t)\\
		\dot \theta_i(t)&=u_i(t)
	\end{split}
\end{align}
where $(x_i,y_i)\in\mathbb{R}^{\mathrm{2}}$ denotes the position of the $i$th pursuer in the frame $X_IOY_I$.

Considering the system of the group of pursuers, the state equations can be obtained as
\begin{equation}
	\label{equation:multi-m}
	\begin{cases}
		\dot x_1(t)=\cos\theta_{M_1}(t)\\
		\dot y_1(t)=\sin\theta_{M_1}(t)\\
		\dot \theta_1(t)=u_1(t)\\
		\dot x_2(t)=\cos\theta_{M_2}(t)\\
		\dot y_2(t)=\sin\theta_{M_2}(t)\\
		\dot \theta_2(t)=u_2(t)\\
		\cdots\\
		\dot x_N(t)=\cos\theta_{M_N}(t)\\
		\dot y_N(t)=\sin\theta_{M_N}(t)\\
		\dot \theta_N(t)=u_N(t)
	\end{cases}.
\end{equation}

\subsection{Optimal Cooperative Intercept Problem}
For a single pursuer, it was proved that the global solution of free-time minimum-effort nonlinear optimal intercept problem does not exist \cite{chen2019nonlinear}. This non-existence of global solutions holds true in the multi-pursuer scenarios. Thus, the objective function with a linear combination of control effort and engagement duration is considered. Then, let us consider the following optimal cooperative intercept problem (OCIP):

\emph{Problem 1(OCIP)}\quad Given an initial condition 
\begin{align}
	(S_1,S_2,\cdots,S_i,\cdots,S_N)
\end{align}
 where $S_i= (x_{i0},y_{i0},\theta_{M_{i0}})\in \mathbb{R}^{\mathrm{2}}\times[-\pi,\pi]$, and a free final time $t_f > 0$, the OCIP consists of steering the system in Eq.~(\ref{equation:multi-m}) by a measurable control vector $(u_1(\cdot),u_2(\cdot),\cdots,u_N(\cdot))$ on $[0,t_f]$ from the initial state $(S_1,S_2,\cdots,S_i,\cdots,S_N)$ to the final point $(0,0)$, with the final heading angle satisfying
\begin{equation}
	\label{equation:constrain}
	\theta_{M_{i+1}}(t_f)-\theta_{M_i}(t_f)=\delta
\end{equation}
such that
\begin{equation}
	\label{sample:equation7}
	J = \int_{0}^{t_f} \kappa+\frac{1}{2}(1-\kappa)\sum_{i=1}^{N}u_i(t)^2dt
\end{equation}
is minimized where $\kappa \in (0,1)$ is a weighting factor.

%%%%%%%%%%%%%%%%%%%%%%%%%%%%%%%%%%%%%%%%%%%%%
\section{Parametrization of Optimal Guidance Law}
%%%%%%%%%%%%%%%%%%%%%%%%%%%%%%%%%%%%%%%%%%%%%
\subsection{Pontryagin’s Maximum Principle}
Denote by $p_{x_i}$, $p_{y_i}$, and $p_{\theta_i}$ as the costate variables of $x_i$, $y_i$, and $\theta_{M_i}$, respectively. Then, the Hamiltonian for the OCIP is expressed as
\begin{equation}
	\begin{split}
		\label{sample:Ham}
		H = &\sum_{i=1}^N p_{x_i}\cos\theta_{M_i}+\sum_{i=1}^N p_{y_i}\sin\theta_{M_i}+\sum_{i=1}^N p_{\theta i}u_i+\\
		&p^0\left[\kappa+\frac{1}{2}(1-\kappa)\sum_{i=1}^N u_i^2\right]
	\end{split}.
\end{equation}
In view of remark 2 from \cite{chen2019nonlinear}, we can consider $p^0=-1$ in the remainder of the paper.

According to PMP, for $t\in[0,t_f]$, it holds that
\begin{align}
	\begin{split}
		\dot p_{x_i} &= -\frac{\partial H}{\partial x_i} = 0\\
		\dot p_{y_i} &=-\frac{\partial H}{\partial y_i} = 0\\
		\dot p_{\theta_i} &= -\frac{\partial H}{\partial \theta_{M_i}} = p_{x_i}\sin\theta_{M_i}-p_{y_i}\cos\theta_{M_i}
	\end{split}
\end{align}
and
\begin{equation}
	\label{equation:Hu}
	\frac{\partial H}{\partial u_i} = 0.
\end{equation}

Explicitly rewriting Eq.~(\ref{equation:Hu}) leads to
\begin{equation}
	\label{equation:u}
	u_i(t)= \frac{p_{\theta_i}(t)}{1-\kappa}.
\end{equation}

Because of the final heading angle in Eq.~(\ref{equation:constrain}), the transversality condition implies
\begin{equation}
	\label{equation:transversality}
	p_{\theta_{N}}(t_f) = -\sum_{i=1}^{N-1}p_{\theta_i}(t_f).
\end{equation}
As the final time is free, along the optimal trajectory, it holds that
\begin{equation}
	\label{equation:Hamiltonian}
	H \equiv 0.
\end{equation}

\subsection{Parameterized System}
Utilizing the canonical equations derived by PMP, we can establish a parameterized system with the gathered states of pursuers as
\begin{align}
	\label{eq:para_sys}
	\begin{split}
		&\begin{cases}
			\dot x_1(t)=\cos\theta_{M_1}(t)\\
			\dot y_1(t)=\sin\theta_{M_1}(t)\\
			\dot \theta_1(t)=u_1(t)\\
			\dot x_2(t)=\cos\theta_{M_2}(t)\\
			\dot y_2(t)=\sin\theta_{M_2}(t)\\
			\dot \theta_2(t)=u_2(t)\\
			\cdots\\
			\dot x_N(t)=\cos\theta_{M_N}(t)\\
			\dot y_N(t)=\sin\theta_{M_N}(t)\\
			\dot \theta_N(t)=u_N(t)
		\end{cases}\\
		&\begin{cases}
			\dot p_{x_1} = 0\\
			\dot p_{y_1} = 0\\
			\dot p_{\theta_1} = p_{x_1}\sin\theta_{M_1}-p_{y_1}\cos\theta_{M_1}\\
			\dot p_{x_2} = 0\\
			\dot p_{y_2} = 0\\
			\dot p_{\theta_i} = p_{x_2}\sin\theta_{M_2}-p_{y_2}\cos\theta_{M_2}\\
			\cdots\\
			\dot p_{x_N} = 0\\
			\dot p_{y_N} = 0\\
			\dot p_{\theta_i} = p_{x_N}\sin\theta_{M_N}-p_{y_N}\cos\theta_{M_N}
		\end{cases}
	\end{split}.
\end{align}

In the following section, a mapping will be established for an FNN to generate the solution of OCIP.
%%%%%%%%%%%%%%%%%%%%%%%%%%%%%%%%%%%%%%%%%%%%%%%%%%%%%
\section{Generation of Dataset}
%%%%%%%%%%%%%%%%%%%%%%%%%%%%%%%%%%%%%%%%%%%%%%%%%%%%%
%\subsection{Backward Propagation}
The parameterized system in Eq.(\ref{eq:para_sys}) can be propagated backward with an initial value satisfying the terminal constraints
\begin{align}
	\begin{split}
		&x_i(t_f)=0\\
		&y_i(t_f)=0\\
		&\theta_{M_{i+1}}(t_f)-\theta_{M_i}(t_f)=\delta
	\end{split}
\end{align}
and transversality conditions
\begin{align}
	\begin{split}
		&p_{\theta_{N}}(t_f) = -\sum_{i=1}^{N-1}p_{\theta_i}(t_f)\\
		&H(t_f) = 0
	\end{split}.
\end{align}
Thus, the nonlinear optimal guidance commands are able to be generated to satisfy the necessary conditions of optimality.

%\subsection{Mapping of Dataset}
Given any appropriate initial value, we can obtain a list of nonlinear optimal guidance commands by solving an initial value problem in Eq.~(\ref{eq:para_sys}). As a result, by sampling some pairs in feasible regions as initial conditions, we are able to use the initial value problem in Eq.~(\ref{eq:para_sys}) to generate sampled data for optimal guidance law. Then, the dataset for the mapping from the flight state to the corresponding optimal feedback control can be immediately obtained. 

\begin{figure}	
	\centering
	\includegraphics[width=3.25in]{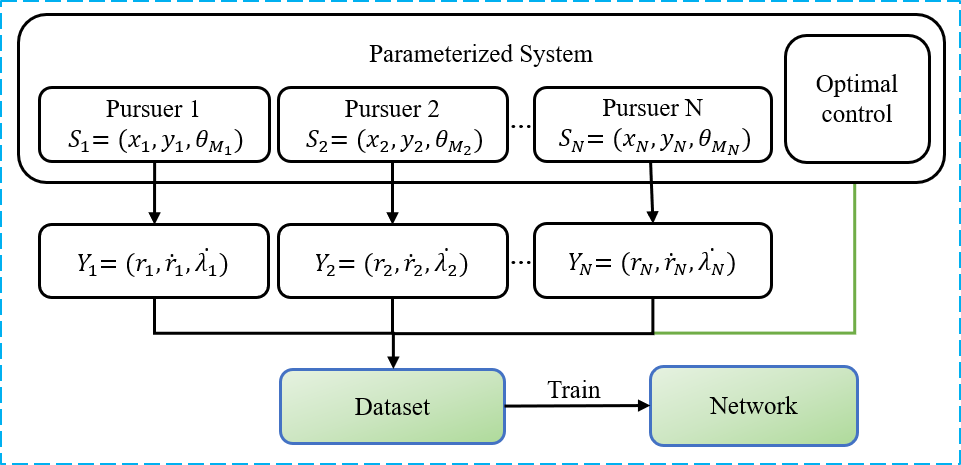}\\
	\caption{\textbf{ Diagram of dataset generation.}}
	\label{fig:train}
\end{figure}

\begin{figure}	
	\centering
	\includegraphics[width=3.25in]{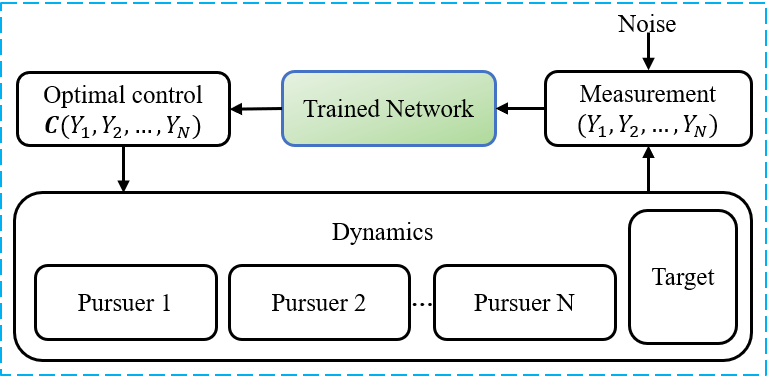}\\
	\caption{\textbf{ Diagram for guidance loop.}}
	\label{fig:guidance}
\end{figure}

In this paper, we assume an ideal condition that the pursuers in the group can communicate with each other. We denote the state of the $i$th pursuer in a polar frame as $Y_i=(r_i,\dot r_i,\dot\lambda_i)$. Let $\C(Y_1,Y_2,\cdots,Y_i,\cdots,Y_N)$ be the optimal feedback control vector of the OCIP at the combined state of pursuers $(Y_1,Y_2,\cdots,Y_i,\cdots,Y_N)$. A simple FNN trained by the dataset is able to approximate the mapping
\begin{equation}
	 (Y_1,Y_2,\cdots,Y_i,\cdots,Y_N) \rightarrow \C(Y_1,Y_2,\cdots,Y_i,\cdots,Y_N)
\end{equation}
as shown by the diagram in Fig.~\ref{fig:train}. Thus, given a combined state of the group of pursuers $(Y_1,Y_2,\cdots,Y_i,\cdots,Y_N)$, each pursuer can use the same FNN to generate the feedback control of the OCIP in a closed-loop guidance system, as shown in Fig.~\ref{fig:guidance}.

%\begin{equation}
%	\label{equation:ode1}
%	\dot r_i(t)=-\cos\sigma_i(t)
%\end{equation}
%\begin{equation}
%	\label{equation:ode2}
%	\dot \sigma_i(t)=\frac{\sin\sigma_i(t)}{r_i(t)}-u_i(t)
%\end{equation}
%\begin{equation}
%	\label{equation:ode3}
%	\dot {\Delta\lambda_i}(t)=\frac{\sin\sigma_{i+1}(t)}{r_{i+1}(t)}-\frac{\sin\sigma_i(t)}{r_i(t)}
%\end{equation}

%$(r_i,\sigma_i,\Delta\lambda_i)$

%\begin{figure*}
%	\centering
%	\includegraphics[width=4in]{Placeholder.eps}
%	\caption{\bf{Here is an example of a figure that spans both columns.}}
%	\label{FlowChart}
%\end{figure*}

%%%%%%%%%%%%%%%%%%%%%%%%%%%%%%%%%%%%%%%%%%%%%%%%%%%%%
\section{Numerical Applications}
%%%%%%%%%%%%%%%%%%%%%%%%%%%%%%%%%%%%%%%%%%%%%%%%%%%%%
This section provides simulations of a two-on-one engagement ($N=2$)  to demonstrate the performance of the proposed nonlinear cooperative optimal guidance strategy. The speed of each pursuer is set as $1000$ m/s. 

Set the weight $\kappa$ as $0.01$ and the relative final heading angle $\delta$ as $10$ deg. We generated about $4.5\times 10^5$ trajectories by backward propagation of the differential equations in Eq.~(\ref{eq:para_sys}). An FNN with three hidden layers (each containing 20 neurons) was trained by the dataset to approximate the optimal feedback guidance commands. Given any feasible input, the trained FNN takes around $0.21$ ms to produce an output on an embedded system with MYC-Y6ULY2 CPU at 528 MHz.

Due to the transversality conditions in Eq.~(\ref{equation:transversality}) and the formula of optimal control in Eq.~(\ref{equation:u}), there must be a one-to-many mapping at the terminal of the OCIP. According to the universal approximation theorem (see, e.g., \cite{hornik1989multilayer,cybenko1989approximation,hornik1991approximation}), the FNN is unable to approximate the mapping. Therefore, we will switch to proportional navigation (PN) once the range between the $i$th pursuer and the target is less than $200$ m.

\begin{figure}	
	\centering
	\includegraphics[width=3.25in]{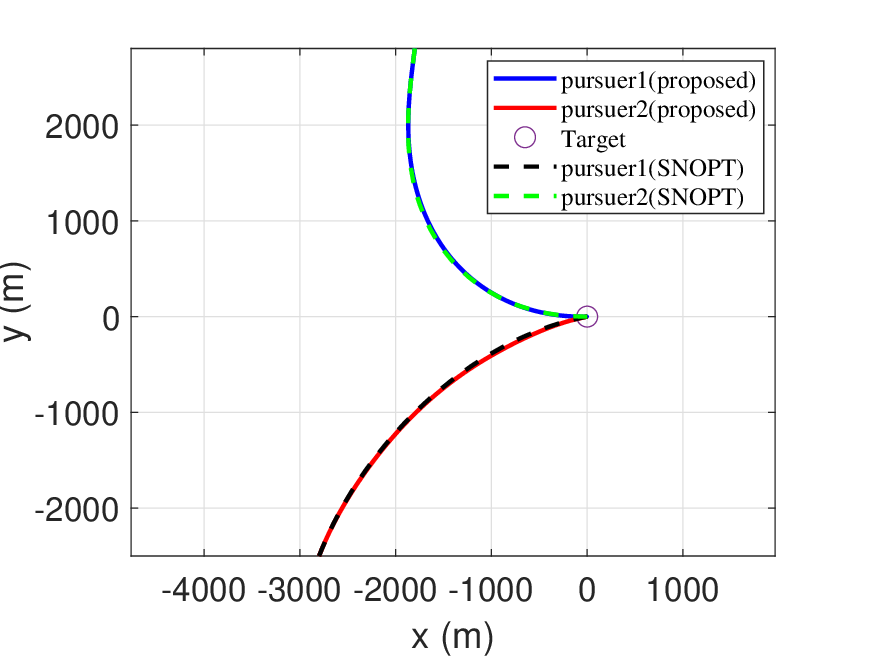}\\
	\caption{\textbf{ Trajectories for case 1.}}
	\label{traj}
\end{figure}

\begin{figure}
	\centering
	\includegraphics[width=3.25in]{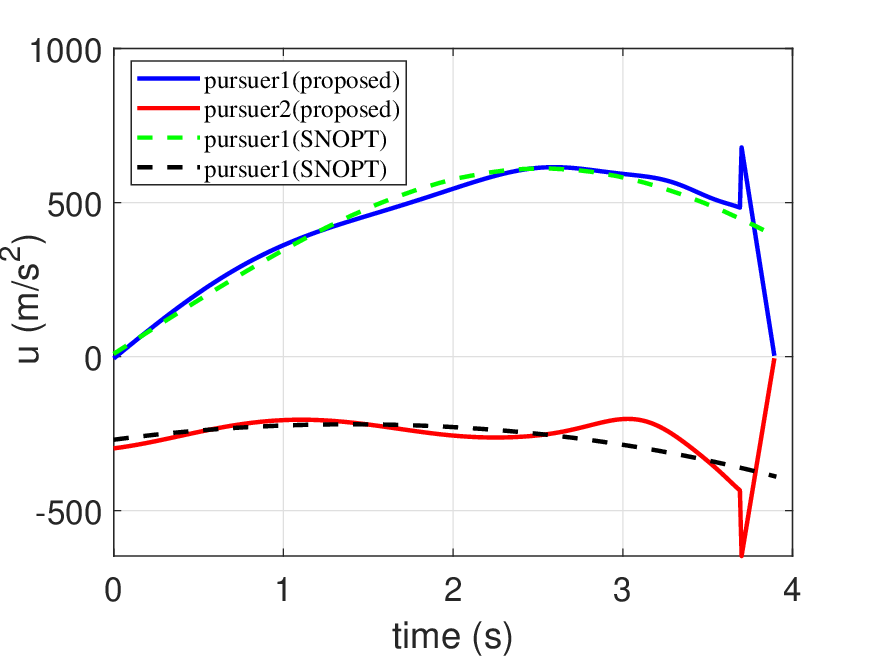}\\
	\caption{\textbf{ Control profiles for case 1.}}
	\label{control}
\end{figure}

\begin{figure}
	\centering
	\includegraphics[width=3.25in]{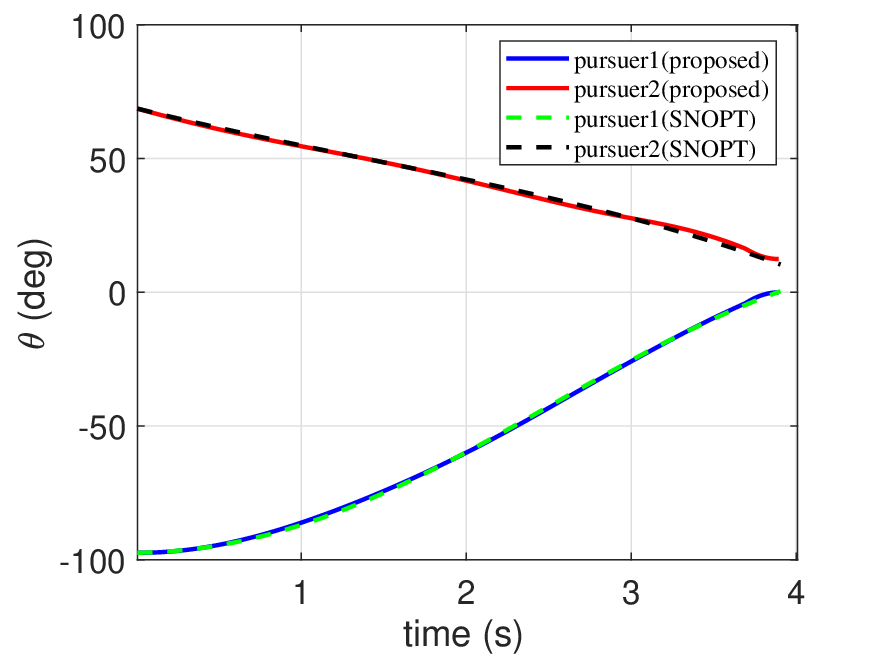}\\
	\caption{\textbf{ Time histories of heading angle for case 1.}}
	\label{heading}
\end{figure}

\begin{figure}
	\centering
	\includegraphics[width=3.25in]{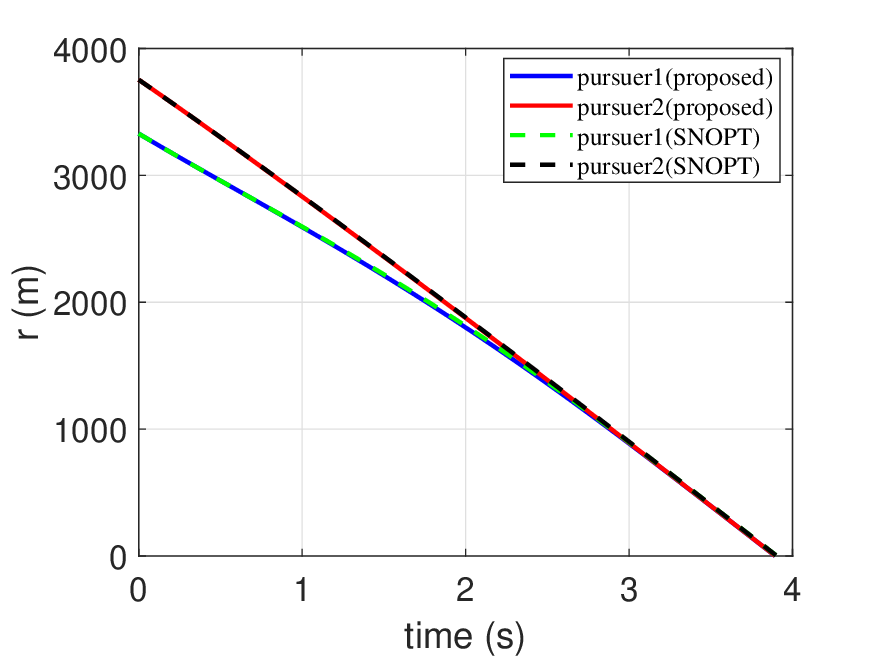}\\
	\caption{\textbf{ Time histories of range for case 1.}}
	\label{range}
\end{figure}

\subsection{Case $1$}
In order to show the capability of the trained network to generate optimal trajectories, an example of proximity is simulated by comparing it with the optimization methods (here, we employed the optimization toolbox: SNOPT). Their initial conditions of them are given in Table \ref{Init1}. The trajectories of the two pursuers are presented in Fig.~\ref{traj}. Their control profiles are demonstrated in Fig.~\ref{control}. The heading angle profiles are presented in Fig.~\ref{heading}. The time histories of range between $i$th pursuer and the target are depicted in Fig.~\ref{range}. Note that the proposed guidance strategy should be switched to PN once the range $r<200$ m. Thus, the control profiles are not continuous for the proposed guidance strategy. It is seen that the trajectories generated by the trained network coincide with those obtained from optimization methods. However, as optimization methods suffer the convergence issue, they cannot be guaranteed to calculate the optimal command within a guidance period.

\begin{table}
\renewcommand{\arraystretch}{1.3}
\caption{\bf Initial Condition of Case 1}
\label{Init1}
\centering
\begin{tabular}{|c|c|c|c|}
\hline
\bfseries state & \bfseries $x$ & \bfseries $y$ & \bfseries $\theta$ \\
\hline\hline
pursuer $\#1$ & -1.8km & 2.8km & -97deg \\
pursuer $\#2$ & -2.8km & -2.5km & 69deg\\

\hline
\end{tabular}
\end{table}

\subsection{Case $2$}
\begin{table}
	\renewcommand{\arraystretch}{1.3}
	\caption{\bf Initial Condition of Case 2}
	\label{Init2}
	\centering
	\begin{tabular}{|c|c|c|c|}
		\hline
		\bfseries state & \bfseries $x$ & \bfseries $y$ & \bfseries $\theta$ \\
		\hline\hline
		pursuer $\#1$ & -3.1km & -2.4km & 178deg \\
		pursuer $\#2$ & -0.25km & -5km & 100deg\\
		
		\hline
	\end{tabular}
\end{table}

\begin{figure}
	\centering
	\includegraphics[width=3.25in]{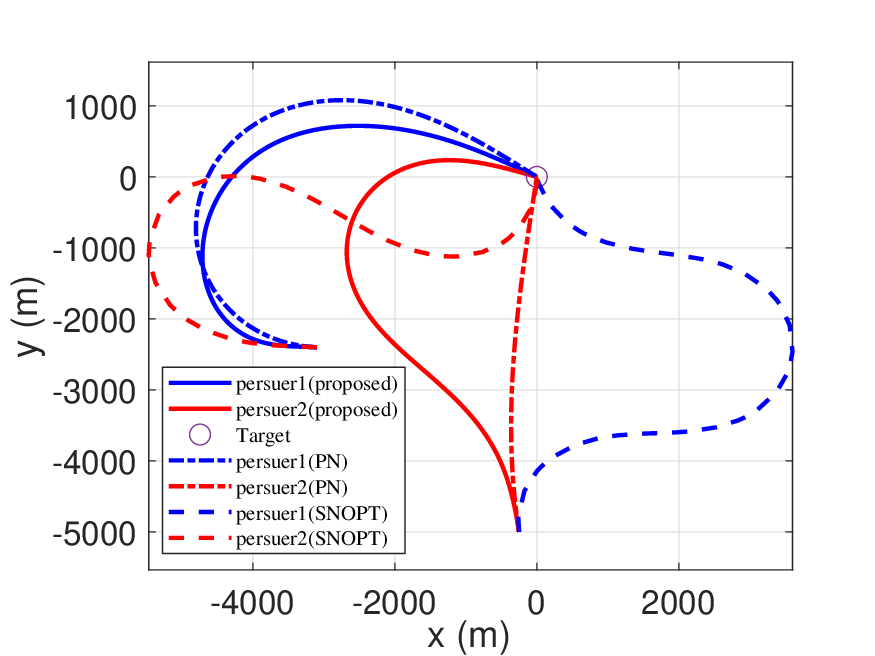}\\
	\caption{\textbf{ Trajectories for case 2.}}
	\label{traj2}
\end{figure}

\begin{figure}
	\centering
	\includegraphics[width=3.25in]{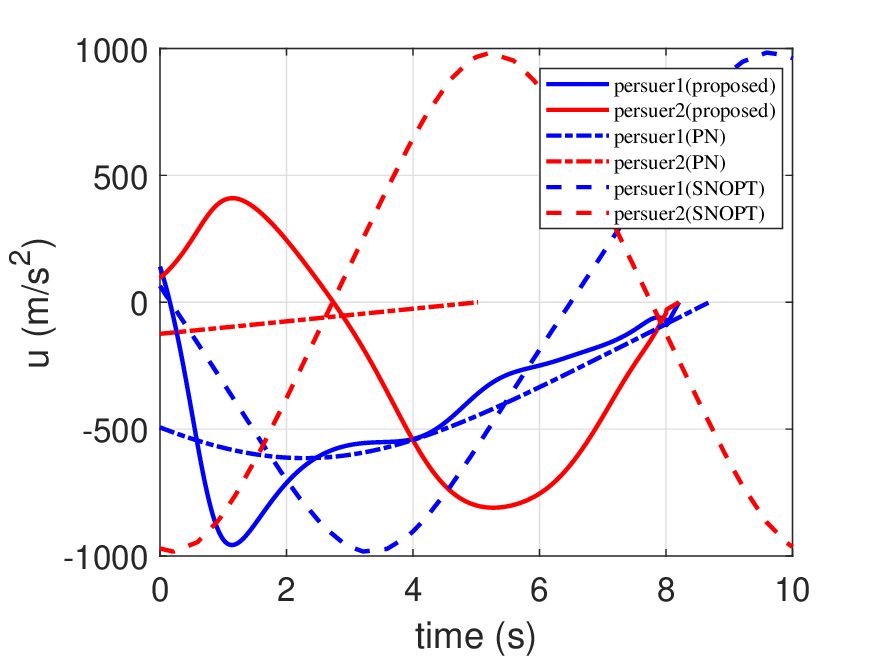}\\
	\caption{\textbf{  Control profiles for case 2.}}
	\label{control2}
\end{figure}

\begin{figure}
	\centering
	\includegraphics[width=3.25in]{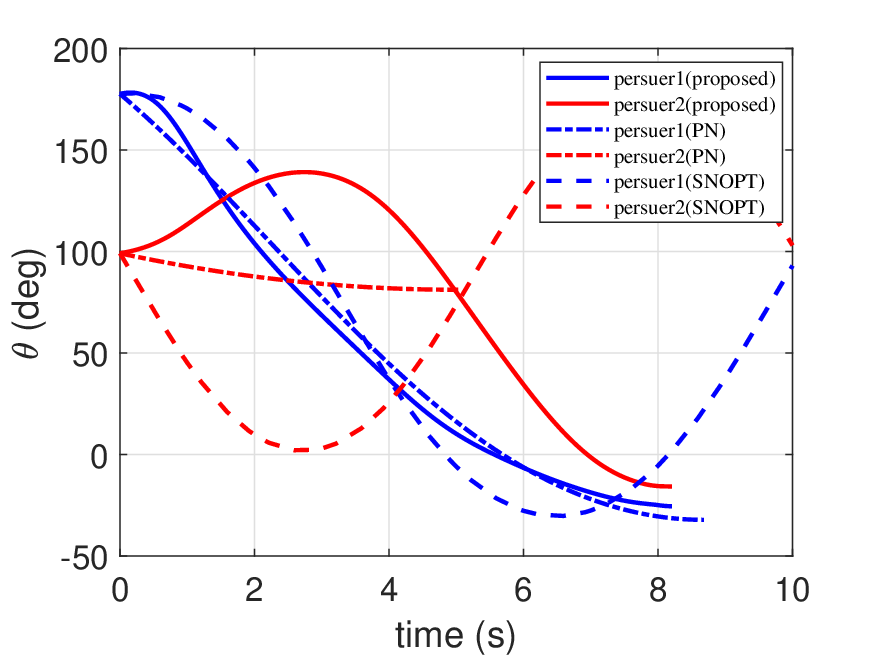}\\
	\caption{\textbf{ Time histories of heading angle for case 2.}}
	\label{heading2}
\end{figure}

\begin{figure}
	\centering
	\includegraphics[width=3.25in]{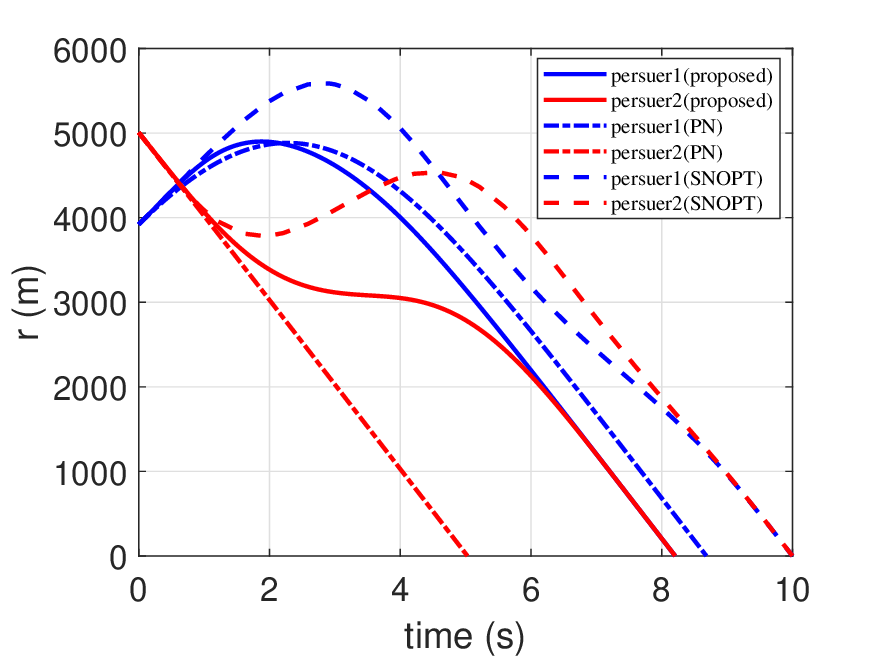}\\
	\caption{\textbf{ Time histories of range for case 2.}}
	\label{range2}
\end{figure}

In this scenario, a simulation of engagement is considered by comparing it with the traditional PN and the optimization methods (SNOPT). Their initial conditions are given in Table \ref{Init2}. The generated trajectories are presented in Fig.~\ref{traj2}, and the corresponding control profiles are reported in Fig.~\ref{control2}. The time histories of heading angle are presented in Fig.~\ref{heading2}. The range profiles are depicted in Fig.~\ref{range2}. It can be seen that in Fig.~\ref{traj2} and Fig.~\ref{heading2}, the constraint of relative intercept angle cannot be satisfied by the pure PN. The optimality of the results is compared with the optimal solutions solved by SNOPT. The total control effort by SNOPT 
\begin{equation}
	\frac{1}{2}\int_{0}^{t_f}u_1^2+u_2^2\ dt=4.3553\times10^6\ m^2/s^3.
\end{equation}
However, the corresponding total control effort by FNN is just $2.0185\times10^6\ m^2/s^3$. It means that the trajectories by FNN are better than that by SNOPT from the perspective of control effort.

%%%%%%%%%%%%%%%%%%%%%%%%%%%%%%%%%%%%%%%%%%%%%%%%%%%%%
\section{Conclusions}
%%%%%%%%%%%%%%%%%%%%%%%%%%%%%%%%%%%%%%%%%%%%%%%%%%%%%
The problem of optimal explicit cooperative guidance for imposing a relative intercept angle was investigated. The method of parametrization for nonlinear optimal guidance law was proposed. Necessary optimality conditions for optimal control policy were established by applying the PMP. Consequently, the optimal trajectories can be propagated from those equations, which can be used to generate the dataset of the mapping from state of the group of pursuers to optimal feedback guidance command. Then, a  simple FNN trained from the dataset was utilized to generate the cooperative guidance command. Finally, two cases of simulations were presented, indicating that the proposed optimal guidance law can generate optimal trajectories within a constant time.

%%%%%%%%%%%%%%%%%%%%%%%%%%%%%%%%%%%%%%%%%%%%%%%%%%%%%%%%%%%%%%%%%%%%%%%%%%%%%%%%%%%%%%%%%%%%%%%%%
%\appendices{}              % note there is no {} to put a title. Each appendix has its own title
%%%%%%%%%%%%%%%%%%%%%%%%%%%%%%%%%%%%%%%%%%%%%%%%%%%%%%%%%%%%%%%%%%%%%%%%%%%%%%%%%%%%%%%%%%%%%%%%%
% For a single appendix, use the \appendix{} keyword and do not use the \section command.

%\section{More Information}        % first appendix
%%%%%%%%%%%%%%%%%%%%%%%%%%
%This is the first appendix. 

%\subsection{Comments}
%If you have only one appendix, use the ``appendix'' keyword.

%\subsection{More Comments}
%Use section and subsection keywords as usual.

%\section{Yet More Information}    % second appendix
%%%%%%%%%%%%%%%%%%%%%%%%%%%%%%
%This is the second appendix.

%%%%%%%%%%%%%%%%%%%%%%%%%%%%%%%%%%%%%%%%%%%%%%%%%%%%%%%%%%%%%%%%%%%%%%%%%%%%%%%%%%%%%%%%%%%%%%%%%%%%%%
\acknowledgments
This work was supported by the National Natural Science Foundation of China (61903331, 62088101), Key Research and Development Program of Zhejiang Province (2020C05001).

%%%%%%%%%%%%%%%%%%%%%%%%%%%%%%%%%%%%%%%%%%%%%%%%%%%%%%%%%%%%%%%%%%%%%%%%%%%%%%%%%%%%%%%%%%%%%%%%%%%%%%
\bibliographystyle{IEEEtran}
%\bibliography{IEEEabr,MyBibFile}

%%%%%%%%%%%%%%%%%%%%%%%%%%%%%%%%%%%%%%%%%%%%%%%%%%%%%%%%%%%%%%%%%%%%%%%%%%%%%%%%%%%%%%%%%%%%%%%%%%%%%%

\end{document}